\documentclass[12pt]{article}
\usepackage[cp1251]{inputenc}
\usepackage{amsmath}
\usepackage{amssymb}%,cite}
\usepackage[english]{babel}
\usepackage[T2A]{fontenc}

\begin{document}

%{\bf 11.02.2020}
\begin{large}

\centerline{\textbf{Degenerate Bianchi trans\-for\-ma\-ti\-ons for
three-dimensional }}

\centerline{\textbf{ pseudo-spherical submanifolds in $\mathbb
R^5$}}

\end{large}

\bigskip

\centerline{A.~A.~Borisenko, V.~O.~Gorkavyy  }

\centerline{(B.~Verkin Institute for Low Temperature Physics and
Engineering, Kharkiv, Ukraine)}
\bigskip

{\small {\bf Abstract.} Three-dimensional pseudo-spherical
submanifolds in $\mathbb R^5$, whose Bianchi
trans\-for\-ma\-ti\-ons are degenerate of rank 2, are studied. A
complete description of such submanifolds is obtained in the case
where the Bianchi trans\-for\-ma\-ti\-ons are holonomically
degenerate.}

\bigskip

\textbf{Introduction}

\medskip

A regular $n$-dimensional submanifold $F^n$ in $(n+p)$-dimensional
Euclidean space $\mathbb R^{n+p}$ is called pseudo-spherical if
its sectional curvature is constant negative, $K\equiv -1$. Any
such submanifold represents a domain of the hyperbolic space $H^n$
isometrically immersed into $\mathbb R^{n+p}$. Consequently, it
can be parameterized by horospherical coordinates $(u_1, ...,
u_n)$ so that its metric form is
\begin{equation}\label{metric}
ds^2 = du_1^2 + e^{-2u_1}\left(u_2^2+...u_n^2\right).
\end{equation}
Let $x(u_1,..., u_n)$ denote the position vector of $F^n$ in $\mathbb R^{n+p}$.

The Bianchi trans\-for\-ma\-ti\-on of $F^n$ corresponding to the chosen horospherical coordinates
is defined by the formula
\begin{equation}\label{def}\tilde x = x +\frac{\partial x}{\partial u_1}.
\end{equation}

The classical theory of pseudo-spherical submanifolds and their Bianchi trans\-for\-ma\-ti\-ons
deals with the case of $n$-dimensional submanifolds in $\mathbb R^{2n-1}$. In this situation, the
Bianchi trans\-for\-ma\-ti\-on has remarkable properties. For instance, if the resulting
vector-function $\tilde x(u_1, ..., u_n)$ represents a regular $n$-dimensional submanifold, then
this submanifold is pseudo-spherical and its sectional curvature is $\tilde K\equiv -1$, see
\cite{Bianchi}, \cite{Am1}, \cite{M}, \cite{T}. Thus, the concept of Bianchi trans\-for\-ma\-ti\-on
provides us with an effective way for constructing new pseudo-spherical submanifolds from a given
one.

Moreover, the theory of $n$-dimensional pseudo-spherical
submanifolds in $\mathbb R^{2n-1}$ has a remarkable analytical
interpretation in the frames of the theory of solitons, where
these submanifolds are represented by solutions of a specific
integrable system of non-linear pde's generalizing the celebrated
sine-Gordon equation, see \cite{T}, \cite{TT1}, \cite{TT2}.

As for the case of $p>n-1$, it quite differs from the classical
one where $p=n-1$. For instance, in this non-classical case the
pseudo-sphericity is not obliged to be preserved under the Bianchi
trans\-for\-ma\-ti\-on, see  \cite{AmS}, \cite{G1}, \cite{G2}.

Generically, the vector-function $\tilde x(u_1,..., u_n)$ defined
by (\ref{def}) represents an $n$-dimensional sub\-ma\-ni\-fold in
$\mathbb R^{n+p}$, which is regular almost everywhere. However, in
some particular cases the Bianchi trans\-for\-ma\-ti\-on
degenerates in the sense that $\tilde x(u_1,..., u_n)$ describes a
submanifold whose dimension is less than $n$.

The simplest example of a degenerate Bianchi trans\-for\-ma\-ti\-on is provided by the famous
pseudo-sphere (Beltrami surface), which is a surface of revolution in $\mathbb R^3$ obtained by
rotating a tractrix. This surface admits a degenerate Bianchi trans\-for\-ma\-ti\-on to a straight
line.

A complete description of $n$-dimensional pseudo-spherical submanifolds in $\mathbb R^{n+p}$
admitting degenerate Bianchi trans\-for\-ma\-ti\-ons to one-dimensional curves was proved recently
in \cite{GN1}, see also \cite{BMAG}, \cite{GN2}. These submanifolds, which are called generalized
Beltrami surfaces, are submanifolds of revolution obtained by rotating particular curves,
generalized tractrices.

The aim of this paper is to analyze the degenerate Bianchi
trans\-for\-ma\-ti\-ons of rank 2 which by definition transform
$n$-dimensional submanifolds with $n>2$ to two-dimensional
surfaces. In this paper we consider the case of three-dimensional
pseudo-spherical submanifolds in $\mathbb R^5$.

So, let $F^3$ be a pseudo-spherical submanifold in $\mathbb R^5$. Suppose that $F^3$ admits a
Bianchi trans\-for\-ma\-ti\-on which degenerates and transforms $F^3$ into a two-dimensional
surface, $\tilde F^2$. What can one say about the initial submanifold $F^3$ as well as about the
transformed surface $\tilde F^2$?

Evidently, the kernel of (the differential of)  the degenerate Bianchi trans\-for\-ma\-ti\-on in
question represents a well-defined one-dimensional distribution in the tangent bundle $TF^3$. It
turns out that this distribution is tangent to the coordinate horospheres $u_1=const$, and its
integral trajectories called the null curves of the degenerate Bianchi trans\-for\-ma\-ti\-on
belong to the coordinate horospheres $u_1=const$.

Thus, due to the degeneracy of the Bianchi transformation, there
are three well-defined one-dimensional distributions in $TF^3$:
the distribution spanned by $\frac{\partial x}{\partial u_1}$, the
kernel of the Bianchi transformation, and a one-dimensional
distribution orthogonal to the two distributions above.

\medskip

{\bf Definition.} The degenerate Bianchi transformation is referred to as {\it holonomically
degenerate} if the three one-dimensional distributions in question, as well as the three
two-dimensional distributions determined by them in $TF^3$, are integrable in the sense of the
classical Frobenius theorems, i.e., the Lie brackets of vector fields $\xi_1$, $\xi_2$, $\xi_3$
spanning these three one-dimensional distributions satisfy $\left[ \xi_i, \xi_j \right] \in span (
\xi_i, \xi_j )$, $1\leq i < j\leq 3$.

\medskip

The main result of the paper provides a complete description of three-dimensional
pseudo-sphe\-ri\-cal submanifolds in $\mathbb R^5$ which admit holonomically degenerate Bianchi
trans\-for\-ma\-ti\-ons of rank 2.

\medskip

{\bf Theorem.} {\it Let $F^3$ be a pseudo-spherical submanifold in
$\mathbb R^5$. Suppose that $F^3$ admits a ho\-lo\-no\-mi\-cal\-ly
degenerate Bianchi trans\-for\-ma\-ti\-on of rank 2, $F^3\to\tilde
F^2$.

Then the resulting surface $\tilde F^2$ belongs to a
three-dimensional subspace $\mathbb R^3\subset \mathbb R^5$ and
has constant negative Gauss curvature, $\tilde K\leq -1$.

Moreover, the following description holds true:

1) If $\tilde K\equiv -1$, then $F^3$ is represented by a position vector
\begin{equation*}
x(u_1, u_2, u_3) = \left( \bar x(u_1,u_2), e^{-u_1}\cos u_3,
e^{-u_1}\sin u_3 \right),
\end{equation*}
where $\bar x(u_1,u_2)$ is a generical vector-function
representing a two-dimensional surface in $\mathbb R^3$ with the
first fundamental form
\begin{equation*}
d\bar x^2 = (1-e^{-2u_1})du_1^2+e^{-2u_1}du_2^2.
\end{equation*}
The null-curves of the Bianchi trans\-for\-ma\-ti\-on in question are parallel straight lines
$u_2=const$ on the coordinate horospheres $u_1=const$.

2) If $\tilde K< -1$, then $F^3$ is represented by a position vector
\begin{equation*}
x(v_1, v_2, v_3) = \left( \bar x(v_1,v_2), \frac{1}{a}e^{-v_1}v_2\cos a v_3,
\frac{1}{a}e^{-v_1}v_2\sin a v_3 \right),
\end{equation*}
where $\bar x(v_1,v_2)$ is a generical vector-function
representing a two-dimensional surface in $\mathbb R^3$ with the
first fundamental form
\begin{equation*} d\bar x^2 =
(1-e^{-2v_1}\frac{v_2^2}{a^2})dv_1^2+
2e^{-2v_1}\frac{v_2}{a^2}dv_1dv_2 +
e^{-2v_1}(1-\frac{1}{a^2})dv_2^2,
\end{equation*} the constant $a>1$ is
well-defined from the relation $\tilde K= -\frac{a^2}{a^2-1}$, and
the coordinates $(v_1, v_2, v_3)$ are related to the horospherical
coordinates $(u_1, u_2, u_3)$ by $u_1=v_1$, $u_2=v_2\cos v_3$,
$u_3=v_2\sin v_3$. The null-curves of the Bianchi
trans\-for\-ma\-ti\-on in question are concentric circles
$v_2=const$ on the coordinate horospheres $u_1=const$.}

\medskip

For the moment we don't know what happens if we remove the assumption of holonomicity. Most likely,
in the general case the resulting surface $\tilde F^2$ will no longer be pseudo-spherical. However
we believe that the following conjecture holds true: if a three-dimensional pseudo-spherical
submanifold $F^3\subset \mathbb R^5$ admits a degenerate Bianchi transformation to a
two-dimensional surface $\tilde F^2$ so that $\tilde F^2$ has constant negative Gauss curvature and
belongs to a three-dimensional subspace $\mathbb R^3\subset\mathbb R^5$, then the degenerate
Bianchi transformation in question is holonomically degenerate.

\bigskip

\bigskip

\textbf{1. General degenerate Bianchi trans\-for\-ma\-ti\-ons of rank two }

\medskip

Let $F^3$ be a three-dimensional pseudo-spherical submanifold in the five-dimensional Euclidean
space $\mathbb R^5$. The pseudo-sphericity means that the sectional curvature of $F^3$ is constant
negative, $K\equiv -1$. We suppose that $F^3\subset \mathbb R^5$ is represented locally by a
position vector $x(u_1,u_2,u_3)$ in terms of {\it horospherical} coordinates $(u_1,u_2,u_3)$ so
that the first (metric) fundamental form of $F^3$ is given by (\ref{metric}). The coordinate system
in question is orthogonal, the $u_1$-coordinate curves represent a family of geodesics in $F^3$,
whereas the coordinate surfaces $u_1=const$ are (domains on) horospheres. Each of these (domains
on) horospheres is isometric to (a domain of) the Euclidean plane, and the coordinates $(u_2, u_3)$
are Cartesian.

The horospherical coordinates in $F^3$ being fixed, consider the corresponding Bianchi
trans\-for\-ma\-ti\-on of $F^3$ that is defined by the formula (\ref{def}). Notice that different
horospherical coordinate systems on $F^3$ generate different Bianchi trans\-for\-ma\-ti\-ons.

Generically, the resulting vector-function $\tilde x$ in (\ref{def}) depends on the all three
arguments $u_1$, $u_2$, $u_3$ and therefore represents a three-dimensional submanifold in $\mathbb
R^5$.

However, we are interested in the particular case where $\tilde x
(u_1, u_2, u_3)$ represents not a three-dimensional
sub\-ma\-ni\-fold, but a two-dimensional surface.

More precisely, we say that the Bianchi trans\-for\-ma\-ti\-on is {\it degenerate of rank 2} and
hence produces a two-dimensional surface if
\begin{equation}\label{def2}
\dim\hbox{span}\left\{ \frac{\partial\tilde x}{\partial u_1}, \frac{\partial\tilde
x}{\partial u_2}, \frac{\partial\tilde x}{\partial
u_3}\right\}\equiv 2.\end{equation}

In order to differentiate $\tilde x(u_1, u_2, u_3)$ given by
(\ref{def}) and then verify (\ref{def2}), we need to use the
classical Gauss-Weingarten equations, see \cite{Am}, \cite{Eis}:
\begin{equation*}\label{p1_3_1}
\frac{\partial^2 x}{\partial u_i \partial u_j} =
\sum\limits_k\Gamma_{ij}^k \frac{\partial x}{\partial u_k} +
\sum\limits_\sigma  b^\sigma_{ij} n_\sigma,\quad\quad 1\leq i,j
\leq 3,
\end{equation*}
\begin{equation*}\label{p1_3_2}
\frac{\partial n_\sigma}{\partial u_i } = - \sum\limits_{j,
k}b^\sigma_{ij}g^{jk}\frac{\partial x}{\partial u_k} +
\mu_{\sigma\nu\vert i}n_\nu,\quad\quad 1\leq i \leq 3, \quad\quad
1\leq \sigma\not= \nu \leq 2,
\end{equation*}
where $n_1$, $n_2$ are orthogonal unit vector fields normal to $F^3$, $g_{ij}$ and $\Gamma_{ij}^k$
denote the coefficients of the first fundamental form (\ref{metric}) and the Christoffel symbols of
$F^3$, $b_{ij}^\sigma$ stand for the coefficients of the second fundamental form of $F^3$ with
respect to $n_\sigma$, and $\mu_{\sigma\nu\vert i}$ are the torsion coefficients of the normal
frame $n_1$, $n_2$. It is easy to see that the non-zero Christoffel symbols of the metric form
(\ref{metric}) are $\Gamma_{12}^2=\Gamma_{13}^3=-1$, $\Gamma_{22}^1=\Gamma_{33}^1=e^{-2u_1}$.

By differentiating (\ref{def}) with the help of Gauss-Weingarten equations, we obtain
\begin{equation}\label{p1_4_1}
\frac{\partial\tilde x}{\partial u_1} = \frac{\partial x}{\partial
u_1}  + b_{11}^1 n_1 + b_{11}^2 n_2 ,
\end{equation}
\begin{equation}\label{p1_4_2}
\frac{\partial\tilde x}{\partial u_2} = \quad\quad\quad b_{12}^1
n_1 + b_{12}^2 n_2 ,
\end{equation}
\begin{equation}\label{p1_4_3}
\frac{\partial\tilde x}{\partial u_3} = \quad\quad\quad b_{13}^1
n_1 + b_{13}^2 n_2.
\end{equation}

Therefore, the defining condition (\ref{def2}) rewrites as
follows:
\begin{equation}\label{p1_5}
rank\left(  \begin{array}{cc} b_{12}^1& b_{12}^2\\ b_{13}^1& b_{13}^2\\ \end{array} \right) \equiv
1.
\end{equation}

Clearly, the relation (\ref{p1_5}) means that the normal vector fields $b_{12}^1n_1+b_{12}^2n_2$
and $b_{13}^1n_1+b_{13}^2n_2$ are mutually collinear and determine a well-defined one-dimensional
distribution in the normal bundle $NF^3$. Without loss of generality, we will specify the choice of
the orthonormal frame $n_1$, $n_2$ on $F^3$ so that the distribution above is directed along $n_2$,
and hence we have
\begin{equation}\label{p1_5m1}
b_{12}^1\equiv 0, \quad b_{13}^1\equiv 0,
\end{equation}
\begin{equation}\label{p1_5m2}
(b_{12}^2)^2+(b_{13}^2)\not= 0
\end{equation}
Thus, by specifying the normal vectors, we may replace
(\ref{p1_5}) by (\ref{p1_5m1})-(\ref{p1_5m2}).

From the geometric point of view, the relations (\ref{p1_5m1})-(\ref{p1_5m2}) mean that at every
point in $F^3$ the tangent vector $\frac{\partial x}{\partial u_1}$ represents a principal
direction of the second fundamental form of $F^3$ with respect to the normal vector $n_1$, whereas
it does not represent principal directions of the second fundamental form of $F^3$ with respect to
$n_2$. Notice that this geometric property can be used as an alternative way to define
tree-dimensional pseudo-spherical submanifolds in $\mathbb R^5$ which admit Bianchi
trans\-for\-ma\-ti\-ons degenerate of rank 2.

\bigskip

\textbf{2. Holonomically degenerate Bianchi trans\-for\-ma\-ti\-ons of rank two}

\medskip

In view of (\ref{p1_5m1}), we rewrite
(\ref{p1_4_2})-(\ref{p1_4_3}) as follows:
\begin{equation*}\label{p1_4_2m}
\frac{\partial\tilde x}{\partial u_2} = b_{12}^2 n_2,
\end{equation*}
\begin{equation*}\label{p1_4_3m}
\frac{\partial\tilde x}{\partial u_3} = b_{13}^2 n_2.
\end{equation*}
Hence, we have:
\begin{equation*}\label{p1_8}
- b_{13}^2 \frac{\partial\tilde x}{\partial u_2} +
b_{12}^2\frac{\partial\tilde x}{\partial u_3} = 0.
\end{equation*}

This means that the tangent vector field $- b_{13}^2\frac{\partial x}{\partial u_2}+
b_{12}^2\frac{\partial x}{\partial u_3}$ on $F^3$ spans a well-defined one-dimensional distribution
in the tangent bundle $TF^3$ which vanishes under the Bianchi trans\-for\-ma\-ti\-on. The
distri\-bu\-ti\-on can be treated as {\it the kernel} of the Bianchi trans\-for\-ma\-ti\-on. Its
integral trajectories in $F^3$ are called {\it the null-curves} of the Bianchi
trans\-for\-ma\-ti\-on.

Clearly, the vector field $- b_{13}^2\frac{\partial x}{\partial
u_2}+ b_{12}^2\frac{\partial x}{\partial u_3}$ is tangent to the
coordinate horospheres $u_1=const$ and, consequently its integral
trajectories, the null-curves of the degenerate Bianchi
trans\-for\-ma\-ti\-on, belong to the coordinate horospheres
$u_1=const$.

It is natural to consider the situation where both the null curves  and the $u_1$-coordinate curves
can be incorporated into a three-orthogonal coordinate system on $F^3$.

Namely, together with the vector fields $\xi_1=\frac{\partial x}{\partial u_1}$ and
$\xi_3=-b_{13}^2\frac{\partial x}{\partial u_2}+ b_{12}^2\frac{\partial x}{\partial u_3}$, consider
the third vector field $\xi_2=b_{12}^2\frac{\partial x}{\partial u_2}+ b_{13}^2\frac{\partial
x}{\partial u_3}$ on $F^3$. Taking into account (\ref{metric}), it is easy to see that $\xi_1$,
$\xi_2$ and $\xi_3$ are mutually orthogonal.

By definition, the Bianchi trans\-for\-ma\-ti\-on is {\it
holonomically degenerate}, if the Lie brackets of $\xi_1$, $\xi_2$
and $\xi_3$ satisfy
\begin{equation*}\label{p1_7}
\left[ \xi_i, \xi_j \right] \in span ( \xi_i, \xi_j ),\quad 1\leq i < j\leq 3,
\end{equation*}
and hence there exists locally an orthogonal coordinate system,
$(v_1, v_2, v_3)$, on $F^3$ such that $\xi_1$, $\xi_2$ and $\xi_3$
are tangent to the coordinate curves. By calculating the Lie
brackets in question, it is easy to verify that the following
statement holds true.

\medskip

{\bf Proposition 2.1.} {\it The Bianchi trans\-for\-ma\-ti\-on is
holonomically degenerate if and only if
\begin{equation}\label{p1_9}
\frac{\partial b_{12}^2}{\partial u_1}\, b_{13}^2 - \frac{\partial
b_{13}^2}{\partial u_1}\, b_{12}^2 = 0,
\end{equation}
i.e., if the ratio $\displaystyle{\frac{b^2_{13}}{b^2_{12}}}$ does not depend on $u_1$.}

\medskip

Geometrically, (\ref{p1_9}) means that null-curves situated in different coordinate horospheres
$u_1=const$ are mapped to each other by translations along $u^1$-coordinate curves in $F^3$.

Now let us discuss in more details how the coordinate system $(v_1, v_2, v_3)$ is related to the
original coordinates $(u_1, u_2, u_3)$. By definition of $(v_1, v_2, v_3)$, we have:
\begin{equation}\label{p2_1}
\frac{\partial x}{\partial v_1} = A \frac{\partial x}{\partial
u_1} ,
\end{equation}
\begin{equation}\label{p2_2}
\frac{\partial x}{\partial v_2} = B \left( b_{12}^2\frac{\partial
x}{\partial u_2}+ b_{13}^2\frac{\partial x}{\partial u_3} \right),
\end{equation}
\begin{equation}\label{p2_3}
\frac{\partial x}{\partial v_3} = C \left( -b_{13}^2\frac{\partial
x}{\partial u_2}+ b_{12}^2\frac{\partial x}{\partial u_3} \right),
\end{equation}
where $A$, $B$, $C$ some functions.

These equations rewrite as follows
\begin{equation}\label{p3_1}
\frac{\partial u_1}{\partial v_1} = A, \quad \frac{\partial
u_2}{\partial v_1} = 0, \quad \frac{\partial u_3}{\partial v_1} =
0,
\end{equation}
\begin{equation}\label{p3_2}
\frac{\partial u_1}{\partial v_2} = 0, \quad \frac{\partial
u_2}{\partial v_2} = B b_{12}^2, \quad \frac{\partial
u_3}{\partial v_2} = B b_{13}^2,
\end{equation}
\begin{equation}\label{p3_3}
\frac{\partial u_1}{\partial v_3} = 0, \quad \frac{\partial
u_2}{\partial v_3} = - C b_{13}^2, \quad \frac{\partial
u_3}{\partial v_3} = C b_{12}^2.
\end{equation}

The equations for $u_1(v_1,v_2,v_3)$  in (\ref{p3_1})-(\ref{p3_3})
mean that $A=A(v_1)$. Hence, without loss of generality, we can
set $u_1=v_1$.

Next, it follows from (\ref{p3_1}), that $u_2$ and $u_3$ don't
depend on $v_1$. Moreover, the compatibility of equations for
$u_2(v_2,v_3)$ and $u_3(v_2,v_3)$ in (\ref{p3_1})-(\ref{p3_3})
leads to the following relations for the functions $B$ and $C$:
\begin{equation*}\label{p4_1}
\frac{\partial} {\partial v_1} \left(B b^2_{12}\right)= 0, \quad
\frac{\partial} {\partial v_1} \left(B b^2_{13}\right)= 0;
\end{equation*}
\begin{equation*}\label{p4_2}
\frac{\partial} {\partial v_1} \left(C b^2_{12}\right)= 0, \quad
\frac{\partial} {\partial v_1} \left(C b^2_{13}\right)= 0;
\end{equation*}
\begin{equation*}\label{p4_3}
\frac{\partial} {\partial v_3} \left(B b^2_{12}\right) +
\frac{\partial} {\partial v_2} \left(C b^2_{13}\right)= 0, \quad
\frac{\partial} {\partial v_3} \left(B b^2_{13}\right) -
\frac{\partial} {\partial v_2} \left(C b^2_{12}\right)= 0.
\end{equation*}
The existence of non-vanishing solutions $B$, $C$ to this system
is procured by the relation (\ref{p1_9}).

Let analyze how the fundamental forms of $F^3$ change if we replace $(u_1, u_2, u_3)$ by $(v_1,
v_2, v_3)$.

The first fundamental form (\ref{metric}) of $F^3$ is written in new coordinates ($v_1$, $v_2$,
$v_3$) as follows, due to (\ref{p2_1})-(\ref{p2_3}):
\begin{equation*}\label{p3_m}
ds^2 = dv_1^2 + e^{-2v_1}\left( (b_{12}^2)^2 + (b_{12}^3)^2\right)
\left(B^2dv_2^2 + C^2 dv_3^2\right).
\end{equation*}
It is easy to see in view of (\ref{p3_1})-(\ref{p3_2}) that
neither $\left( (b_{12}^2)^2 + (b_{12}^3)^2\right) B^2$ nor
$\left( (b_{12}^2)^2 + (b_{12}^3)^2\right) C^2$ depend on $v_1$.
Therefore the metric form of $F^3$ reads as follows:
\begin{equation}\label{p3_m1}
ds^2 = dv_1^2 + e^{-2v_1} \left(a_{22}dv_2^2 +
a_{33}dv_3^2\right),
\end{equation}
where $a_{22}=a_{22}(v_2,v_3)$ and $a_{33}=a_{33}(v_2,v_3)$ are
some functions. Clearly, $a_{22}dv_2^2 + a_{33}dv_3^2$ is the
metric form of the coordinate horospheres $v_1=const$. It is
Euclidean and this impose some additional constraints on
$a_{22}(v_2,v_3)$ and $a_{33}(v_2,v_3)$.

As for the coefficients of the second fundamental forms of $F^3$
in new coordinates $(v_1, v_2, v_3)$, which will be denoted by
${\text b}^\sigma_{ij}$, it is easy to verify that they still
satisfy conditions similar to (\ref{p1_5m1}),
\begin{equation}\label{p3_71}
{\text b}_{12}^1\equiv 0, \quad {\text b}_{13}^1\equiv 0.
\end{equation}
Besides, we get
\begin{equation}\label{p3_72}
{\text b}_{13}^2\equiv 0
\end{equation}
due to (\ref{p2_1})-(\ref{p2_3}), whereas
\begin{equation}\label{p3_73}
{\text b}_{12}^2\not= 0
\end{equation}
in view of (\ref{p1_5m2}).

Moreover, recall that $n$-dimensional pseudo-spherical submanifolds in ($2n-1$)-dimensional
Euclidean space have flat normal connection, for $n=3$ as well as for any other $n\geq 2$, see
\cite{Am}, \cite{B1}-\cite{B2}, \cite{T}. Since $F^3\subset \mathbb R^5$ has flat normal
connection, then the torsion coefficients $\mu_{\sigma\nu\vert j} = \langle \frac{\partial
n_\sigma}{\partial v_j}, n_\nu \rangle$ can be written as follows:
\begin{equation}\label{p3_8}
\mu_{12\vert j} = - \mu_{21\vert j} = \frac{\partial
\varphi}{\partial v_j}, \quad 1\leq j\leq 3,
\end{equation}
where $\varphi=\varphi(v_1,v_2,v_3)$ is some function. Indeed, since the normal connection is flat,
one can equip $F^3$ with an orthonormal normal frame $n_1^*$, $n_2^*$ which is parallel translated
in the normal bundle $NF^3$. The normal frame $n_1$, $n_2$ on $F^3$ specified above is obtained by
rotating $n_1^*$, $n_2^*$, i.e., $n_1 = \cos\varphi\, n_1^* + \sin\varphi \, n_2^*$, $n_2 = -
\sin\varphi\, n_1^* + \cos\varphi\, n_2^*$. Then one can easily get (\ref{p3_8}) by calculating the
torsion coefficients and taking into account that the frame $n_1^*$, $n_2^*$ is parallel translated
in $NF^3$.

\bigskip

\textbf{3. Gauss-Codazzi-Ricci equations for $F^3$ with
holonomically degenerate Bianchi trans\-for\-ma\-ti\-on}

\medskip

The coefficients of the fundamental forms of $F^3\subset \mathbb
R^5$ satisfy the well-known system of Gauss-Codazzi-Ricci
equations \cite{Am}, \cite{Eis}:
\begin{equation}\label{p1_6_1}
R_{ijkl} = {\text b}^1_{ik}{\text b}^1_{jl}-{\text b}^1_{il}{\text
b}^1_{jk} + {\text b}^2_{ik}{\text b}^2_{jl}-{\text
b}^2_{il}{\text b}^2_{jk};
\end{equation}
\begin{equation}\label{p1_6_2}
\frac{\partial {\text b}^\sigma_{ij}}{\partial v_k} +
\sum\limits_{s}\Gamma^s_{ij}{\text b}^\sigma_{sk}+
\sum\limits_{\nu}\mu_{\nu\sigma\vert k}{\text b}^\nu_{ij} =
\frac{\partial {\text b}^\sigma_{ik}}{\partial v_j} +
\sum\limits_{s}\Gamma^s_{ik}{\text b}^\sigma_{sj}+
\sum\limits_{\nu}\mu_{\nu\sigma\vert j}{\text b}^\nu_{ik};
\end{equation}
\begin{equation*}%\label{p1_6_3}
\frac{\partial \mu_{\sigma\nu\vert i}}{\partial v_j} - \sum\limits_{s,l} {\text
b}^\sigma_{is}g^{sl}{\text b}^\nu_{jl}+ \sum\limits_{\alpha} \mu_{\sigma\alpha\vert
i}\mu_{\alpha\nu\vert j} = \frac{\partial \mu_{\sigma\nu\vert j}}{\partial v_i} - \sum\limits_{s,l}
{\text b}^\sigma_{js}g^{sl}{\text b}^\nu_{il}+ \sum\limits_{\alpha} \mu_{\sigma\alpha\vert
j}\mu_{\alpha\nu\vert i},
\end{equation*}
where $R_{ijkl}$ are coefficients of the Riemannian curvature tensor of $F^3$.

Since $F^3\subset \mathbb R^5$ has flat normal connection, then the last group of equations (Ricci
equations) rewrite in a simpler form,
\begin{equation}\label{p1_6_3}
\sum\limits_{s,l} {\text b}^\sigma_{is}g^{sl}{\text b}^\nu_{jl}
 =  \sum\limits_{s,l} {\text b}^\sigma_{js}g^{sl}{\text b}^\nu_{il}.
\end{equation}

Let us discuss how the specific expression (\ref{p3_m1}) for the metric form of $F^3$, the
constraints (\ref{p3_71})-(\ref{p3_73}) on the coefficients of the second fundamental forms of
$F^3$, and the particular representation (\ref{p3_8}) of the torsion coefficients of $F^3$ affect
the Gauss-Codazzi-Ricci equations (\ref{p1_6_1})-(\ref{p1_6_3}). Proceeding step by step, we will
solve these equations and find coefficients of the fundamental forms of $F^3$.

Ricci equation (\ref{p1_6_3}) with $i=1$, $j=3$, $\sigma=1$,
$\nu=2$ is
$$\nonumber
\frac{e^{2v_1}}{a_{22}} {\text b}^1_{32}{\text b}^2_{12}=0,$$
hence we get
\begin{equation*}\label{R45}
{\text b^1_{23}} \equiv 0.
\end{equation*}

Gauss equation (\ref{p1_6_1}) with $i=1$, $j=k=2$, $l=3$ is
$$\nonumber {\text b^2_{12}}{\text b^2_{23}}= 0,$$ hence we
get
\begin{equation*}\label{G8}
{\text b^2_{23}} \equiv 0.
\end{equation*}

Gauss equation (\ref{p1_6_1}) with $i=1$, $j=l=3$, $k=2$ is
$$\nonumber {\text b^2_{12}}{\text
b^2_{33}} = 0,$$ hence we get
\begin{equation*}\label{G33}
{\text b^2_{33}} \equiv 0.
\end{equation*}

Consequently, the $v_3$-coordinate curves, that are the null
curves of the Bianchi transformation, turn out to be lines of
curvatures of $F^3$. Besides, these curves are asymptotic lines on
$F^3$ with respect to $n_2$.

Next, consider Codazzi equation (\ref{p1_6_2}) with $i=3$, $j=1$,
$k=2$, $\sigma=2$, which rewrites as follows:
$$\nonumber  \frac{1}{a_{22}}\, \frac{\partial a_{22}}{\partial v_3} \, {\text b^2_{12}} = 0.
$$
Therefore, we have $\frac{\partial a_{22}}{\partial v_3}=0$, i.e.,
$a_{22}$ depends only on $v_2$. Without loss of generality,  by
applying a scaling $v_2\to \tilde v_2(v_2)$ if necessary, we can
set
\begin{equation*}\label{C15}
a_{22} \equiv 1.
\end{equation*}

Codazzi equation (\ref{p1_6_2}) with $i=1$, $j=2$, $k=3$,
$\sigma=1$ gives us
$$\nonumber {\text b^2_{12}} \frac{\partial\varphi}{\partial v_3} = 0,
$$
hence we obtain
\begin{equation}\label{C34}
\mu_{12\vert 3} = \frac{\partial\varphi}{\partial v_3} \equiv 0.
\end{equation}

Codazzi equation (\ref{p1_6_2}) with $i=1$, $j=3$, $k=2$,
$\sigma=2$ reads
$$\nonumber \frac{\partial{\text b^2_{12}}}{\partial v_3} =
0, $$ hence we have
\begin{equation*}\label{C55}
{\text b^2_{12}} = {\text b^2_{12}} (v_1,v_2)
\end{equation*}

Similarly, from Codazzi equation (\ref{p1_6_2}) with $i=j=1$,
$k=3$, $\sigma=1$ we get
\begin{equation*}\label{C29}
{\text b^1_{11}} = {\text b^1_{11}} (v_1,v_2),
\end{equation*}
from Codazzi equation (\ref{p1_6_2}) with $i=j=1$, $k=3$,
$\sigma=2$ we get
\begin{equation*}\label{C30}
{\text b^2_{11}} = {\text b^2_{11}} (v_1,v_2),
\end{equation*}
from Codazzi equation (\ref{p1_6_2}) with $i=2$, $j=3$, $k=2$,
$\sigma=1$ we get
\begin{equation*}\label{C59}
{\text b^1_{22}} = {\text b^1_{22}} (v_1,v_2),
\end{equation*}
and from Codazzi equation (\ref{p1_6_2}) with $i=2$, $j=3$, $k=2$,
$\sigma=2$ we get
\begin{equation*}\label{C60}
{\text b^2_{22}} = {\text b^2_{22}} (v_1,v_2).
\end{equation*}

Ricci equation (\ref{p1_6_3}) with $i=1$, $j=2$, $\sigma=1$,
$\nu=2$ is rewritten as follows:
$$\nonumber {\text b^2_{12}}\left( {\text
b^1_{22}} e^{2v_1} - {\text b^1_{11}} \right) = 0,$$
hence we have
\begin{equation*}\label{R20}
{\text b^1_{22}} = e^{-2v_1}{\text b^1_{11}}.
\end{equation*}

Gauss equation (\ref{p1_6_1}) with $i=k=1$, $j=l=3$ is
$$\nonumber {\text
b^1_{11}}{\text b^1_{33}} + a_{33} e^{-2v_1} = 0, $$ hence we get
\begin{equation*}\label{G28}
{\text b^1_{33}} = -e^{-2v_1}\frac{a_{33}}{{\text b^1_{11}}},
\end{equation*}
since ${\text b^1_{11}}$ can not vanish.

Now consider Codazzi equation (\ref{p1_6_2}) with $i=j=3$, $k=1$,
$\sigma=1$, which turns out to read as follows:
$$\nonumber \frac{\partial {\text b^1_{11}}}{\partial
v_1} +{\text b^1_{11}} + \left({\text b^1_{11}}\right)^3 = 0.
$$
Solving this differential equation, we obtain:
\begin{equation}\label{C39}
{\text b^1_{11}} = \frac{1}{\sqrt{e^{2v_1} f(v_2) -1}},
\end{equation}
where $f(v_2)$ is a positive function depending only on $v_2$.

In order to determine $f(v_2)$, consider Codazzi equation
(\ref{p1_6_2}) with $i=k=3$, $j=2$, $\sigma=1$, which reduces to
the following equation:
\begin{equation}\label{C68}
a_{33}\frac{df}{dv_2}+\frac{\partial a_{33}}{\partial v_2} f
 = 0.
\end{equation}

At this step of our discussion we have to take into account the coefficient $a_{33}$. Recall that
the metric form $a_{22}dv_2^2+a_{33}dv_3^2$ is Euclidean as the first fundamental form of the
coordinate horospheres $v_1=const$. Since $a_{22}\equiv 1$, its is easy to verify that there are
only two possibilities for $a_{33}$: either $a_{33}=a_{33}(v_3)$ and hence we can set $a_{33}\equiv
1$ by applying a scaling $v_3\to \tilde v_3(v_3)$ if necessary,  or $a_{33}=(v_2 c_1(v_3)+
c(v_3))^2$, where $c_1(v_3)\not= 0$, and hence we can set $a_{33}= (v_2+c(v_3))^2$ by  applying a
scaling $v_3\to \tilde v_3(v_3)$ if necessary. Consider two cases separately.

\medskip

{\it Case 1.} Let $a_{33}\equiv 1$. Then the equation (\ref{C68}) yields $f(v_2)\equiv const$, and
we can set $f(v_2)\equiv 1$ by applying a shift $v_1\to v_1+const$ if necessary.

Next, in this case from the Codazzi equations with $i=j=3$, $k=1$,
$\sigma=2$ and with $i=k=3$, $j=2$, $\sigma=2$, respectively, one
can find
\begin{equation}\label{C48}
\mu_{12\vert 1} \, = \, \frac{\partial\varphi}{\partial v_1} \, = \, \frac{1}{\sqrt{e^{2v_1}-1}} \,
{\text b^2_{11}},
\end{equation}
\begin{equation}\label{C73}
\mu_{12\vert 2} \, = \, \frac{\partial\varphi}{\partial v_2} \, = \, \frac{1}{\sqrt{e^{2v_1}-1}} \,
{\text b^2_{12}}.
\end{equation}

Finally, the remaining equations from the Gauss-Codazzi-Ricci
system reduce to the following three equations:
\begin{equation}\label{C1_1}
{\text b^2_{11}} {\text b^2_{22}} - ({\text b^2_{12}} )^2 = - \frac{1}{e^{2v_1}-1},
\end{equation}
\begin{equation}\label{C1_2}
 \frac{\partial{\text b^2_{11}}}{\partial v_2} - \frac{\partial{\text b^2_{12}}}{\partial v_1}
 + {\text b^2_{12}}\frac{e^{2v_1}}{e^{2v_1}-1} = 0,
\end{equation}
\begin{equation}\label{C1_3}
\frac{\partial{\text b^2_{12}}}{\partial v_2} - \frac{\partial{\text b^2_{22}}}{\partial v_1}
 - {\text b^2_{11}}\frac{1}{e^{2v_1}-1} -{\text b^2_{22}} = 0.\end{equation}

Thus, we obtain that the first fundamental form of $F^3$ is
\begin{equation*}\label{R1}
ds^2 \, = \, dv_1^2+e^{-2v_1} \left(dv_2^2+dv_3^2\right),
\end{equation*}
the second fundamental form of $F^3$ with respect to $n_1$ is
\begin{equation*}\label{R2}
{\text b^1} \, = \, \frac{1}{\sqrt{e^{2v_1}-1}} \, dv_1^2\, + \,\frac{e^{-2v_1}}{\sqrt{e^{2v_1}-1}}
\, dv_2^2 \, - \, e^{-2v_1}\sqrt{e^{2v_1}-1}\, dv_3^2, \end{equation*} the second fundamental form
of $F^3$ with respect to $n_2$ is
\begin{equation*}\label{R3} {\text b^2} \, = \,
{\text b^2_{11}}\, dv_1^2\, + \, 2{\text b^2_{12}}\, dv_1dv_2\, + \, {\text b^2_{22}} \, dv_2^2,
\end{equation*}
 and the torsion coefficients are expressed by
(\ref{C34}), (\ref{C48})-(\ref{C73}),  where ${\text b^2_{11}}(v_1,v_2)$, ${\text
b^2_{12}}(v_1,v_2)$, ${\text b^2_{22}}(v_1,v_2)$ are some functions satisfying the relations
(\ref{C1_1})-(\ref{C1_3}).

\medskip

{\it Case 2.} Let $a_{33}=(v_2+c(v_3))^2$. Then the equation (\ref{C68}) rewrites as follows:
\begin{equation}\label{Feqn}
(v_2+c(v_3))\frac{df}{dv_2} + 2 f = 0.
\end{equation}
If $c(v_3)$ is not constant, then the only solution of (\ref{Feqn}) is $f(v_2)=0$, which
contradicts to the positiveness of $f(v_2)$ required in (\ref{C39}). Hence $c(v_3)$ has to be
constant, $c(v_3)=c_0$. By applying a shift $v_2\to v_2+c_0$ if necessary, we can set $c_0=0$, and
therefore we have $a_{33}=(v_2)^2$. By solving (\ref{Feqn}), it is easy to get
$\displaystyle{f=\frac{f_0}{(v_2)^2}}$, where $f_0$ is a positive constant.

Next, in this case from the Codazzi equations with $i=j=3$, $k=1$,
$\sigma=2$ and with $i=k=3$, $j=2$, $\sigma=2$, respectively, one
can find
\begin{equation}\label{C48_2}
\mu_{12\vert 1} = \frac{\partial\varphi}{\partial v_1} =\frac{v_2 {\text b^2_{11}} - e^{2v_1}{\text
b^2_{12}}}{\sqrt{f_0e^{2v_1}-(v_2)^2}},
\end{equation}
\begin{equation}\label{C73_2}
\mu_{12\vert 2} = \frac{\partial\varphi}{\partial v_2} =\frac{v_2 {\text b^2_{12}} - e^{2v_1}{\text
b^2_{22}}}{\sqrt{f_0e^{2v_1}-(v_2)^2}}.
\end{equation}

Finally, the remaining equations from the Gauss-Codazzi-Ricci
system reduce to the following three equations:
\begin{equation}\label{C2_1}
{\text b^2_{11}} {\text b^2_{22}} - ({\text b^2_{12}} )^2 = - \frac{f_0}{f_0e^{2v_1}-(v_2)^2},
\end{equation}
\begin{equation}\label{C2_2}
 \frac{\partial{\text b^2_{11}}}{\partial v_2} - \frac{\partial{\text b^2_{12}}}{\partial v_1}
 + \left(  f_0{\text b^2_{12}} - v_2{\text
b^2_{22}} \right)\frac{e^{2v_1}}{f_0e^{2v_1}-(v_2)^2} = 0,
\end{equation}
\begin{equation}\label{C2_3}
\frac{\partial{\text b^2_{12}}}{\partial v_2} -
\frac{\partial{\text b^2_{22}}}{\partial v_1}
 + \left(  v_2{\text b^2_{12}} - f_0{\text
b^2_{11}} \right)\frac{1}{f_0e^{2v_1}-(v_2)^2} -{\text b^2_{22}} = 0.\end{equation}

Thus, we obtain that the first fundamental form of $F^3$ is
\begin{equation*}\label{Q1} ds^2 \, = \,
dv_1^2+e^{-2v_1} \left( dv_2^2 \, +\, v_2^2\, dv_3^2\right), \end{equation*} the second fundamental
form of$F^3$ with respect to $n_1$ is
\begin{equation*}\label{Q2}{\text b^1} \, = \, \frac{v_2}{\sqrt{f_0e^{2v_1}-(v_2)^2}}\, dv_1^2
\, + \, e^{-2v_1}\frac{v_2}{\sqrt{f_0e^{2v_1}-(v_2)^2}}\, dv_2^2 \, -\,
e^{-2v_1}v_2\sqrt{f_0e^{2v_1}-(v_2)^2} \, dv_3^2,
\end{equation*} the second fundamental form of $F^3$
with respect to $n_2$ is
\begin{equation*}\label{Q3} {\text b^2} \, = \,
{\text b^2_{11}} \, dv_1^2 \, + \, 2{\text b^2_{12}}\, dv_1dv_2 \, + \, {\text b^2_{22}} \,
dv_2^2,\end{equation*}
 and the torsion coefficients are expressed by (\ref{C34}),
(\ref{C48_2})-(\ref{C73_2}),  where ${\text b^2_{11}}(v_1,v_2)$, ${\text b^2_{12}}(v_1,v_2)$,
${\text b^2_{22}}(v_1,v_2)$ are some functions which satisfy the relations
(\ref{C2_1})-(\ref{C2_3}).

\bigskip

\textbf{4. Example 1}

\medskip

In this section we will construct a particular three-dimensional
pseudo-spherical submanifold in $\mathbb R^5$ whose Bianchi
trans\-for\-ma\-ti\-on is degenerate of rank 2 and leads to a
two-dimensional surface.

Let $\bar F^2$ be a two-dimensional surface in $\mathbb R^3$
represented by a vector-function $\bar x(v_1,v_2)$ so that the
first fundamental form of $\bar F^2$ reads
\begin{equation}\label{ex1_1} d\bar x^2 = (1 - e^{-2v_1})dv_1^2 + e^{-2v_1}dv_2^2.
\end{equation}
Suppose the second fundamental form $\bar b = \bar b_{11}dv_1^2 + 2\bar b_{12}dv_1dv_2 + \bar
b_{22} dv_2^2$ of $\bar F^2$  is not diagonalized,
\begin{equation}\label{ex1_1ss}
\bar b_{12}\not= 0.
\end{equation}
%and hence the tangent vector $\frac{\partial \bar x}{\partial v_1}$  coincides  nowhere with the
%principal directions of $\bar F^2$.

The surface $\bar F^2\subset \mathbb R^3$ being given, consider a three-dimensional submanifold
$F^3$ in $\mathbb R^5$ represented by the vector-function
\begin{equation}\label{ex1_2a}
x(v_1,v_2,v_3) = \left(\bar x(v_1,v_2), e^{-v_1}\cos v_3,  e^{-v_1} \sin v_3 \right).
\end{equation}

Because of (\ref{ex1_1}), the first fundamental form of $F^3$ reads
\begin{equation}\label{ex1_2q}
dx^2 = dv_1^2 + e^{-2v_1}(dv_2^2 + dv_3^2).
\end{equation}
Therefore, the submanifold $F^3$ is pseudo-spherical, its Gauss
curvature is equal to -1, and ($v_1$, $v_2$, $v_3$) are
horospherical coordinates on $F^3$.

\medskip

{\bf Proposition 4.1.} {\it The Bianchi trans\-for\-ma\-ti\-on of $F^3\subset \mathbb R^5$
corresponding to the choice of horosherical coordinates is degenerate and transforms $F^3$ to a
two-dimensional surface $\tilde F^2\subset \mathbb R^5$ with the following properties:

(i) $\tilde F^2$ belongs to a three-dimensional subspace $\mathbb
R^3$ of $\mathbb R^5$;

(ii) $\tilde F^2$ is pseudo-spherical, its Gauss curvature is equal to -1. }

\medskip

\texttt{Proof.} Applying the Bianchi trans\-for\-ma\-ti\-on
(\ref{def}) to the vector-function (\ref{ex1_2a}), it is easy to
get
\begin{equation}\label{ex1_2}
\tilde x = \left(\bar x+ \frac{\partial \bar x}{\partial v_1}, 0,  0 \right).
\end{equation}
Clearly, the vector-function $\tilde x$ depends only on $v_1$, $v_2$
 and takes values in the subspace $\mathbb R^3$ of $\mathbb
R^5$ determined by $x_4=0$, $x_5=0$.

Differentiate (\ref{ex1_2}) and apply the Weingarten equations corresponding to $\bar F^2$. Then we
have:
\begin{equation}\label{ex1_3a} \frac{\partial \tilde x}{\partial v_1} =  (1+\bar\Gamma^1_{11})
\frac{\partial \bar x}{\partial v_1} + \bar\Gamma^2_{11} \frac{\partial \bar x}{\partial v_2} +
\bar b_{11} \bar n = \frac{1}{1 - e^{-2v_1}}\frac{\partial \bar x}{\partial v_1} + \bar b_{11} \bar
n,
\end{equation}
\begin{equation}\label{ex1_3b}
\frac{\partial \tilde x}{\partial v_2} = \bar\Gamma^1_{12} \frac{\partial \bar x}{\partial v_1} +
(1+\bar\Gamma^2_{12}) \frac{\partial \bar x}{\partial v_2} + \bar b_{12} \bar n = \bar b_{12} \bar
n,
\end{equation}
where $\bar\Gamma^1_{11}=\frac{e^{-2v_1}}{1 - e^{-2v_1}}$, $\bar\Gamma^2_{11}=0$,
$\bar\Gamma^1_{12}=0$, $\bar\Gamma^2_{12}=-1$ are the Christoffel symbols and $\bar n$ is the unit
normal of $\bar F^2\subset \mathbb R^3$. Taking into account (\ref{ex1_1ss}), we see that
$\displaystyle{\frac{\partial \tilde x}{\partial v_2}}$ and $\displaystyle{\frac{\partial \tilde
x}{\partial v_2}}$ are non-collinear.

Thus, $\tilde x(v_1,v_2)$ given by (\ref{ex1_2}) represents a two-dimensional surface, $\tilde
F^2$, in the three-dimensional Euclidean space $\mathbb R^3$.

To finish the proof, let us find the Gauss curvature $\tilde K$ of $\tilde F^2$. It is easy to
verify in view of (\ref{ex1_3a})-(\ref{ex1_3b}) that $\displaystyle{\frac{\partial \bar x}{\partial
v_2}}$ is orthogonal to $\displaystyle{\frac{\partial \tilde x}{\partial v_1}}$ and
$\displaystyle{\frac{\partial \tilde x}{\partial v_2}}$ because of (\ref{ex1_1}). Hence the unit
normal to $\tilde F^2\subset \mathbb R^3$ is written as follows:
\begin{equation}\label{ex1_4}
\tilde n = e^{v_1} \frac{\partial \bar x}{\partial v_2}.
\end{equation}
Differentiate (\ref{ex1_4}) with the help of Weingarten equations for $\bar F^2$, and  express
$\displaystyle{\frac{\partial \bar x}{\partial v_1}}$, $\displaystyle{\frac{\partial \bar
x}{\partial v_2}}$, $\bar n$ in terms of $\displaystyle{\frac{\partial \tilde x}{\partial v_1}}$,
$\displaystyle{\frac{\partial \tilde x}{\partial v_2}}$, $\tilde n$ by using
(\ref{ex1_3a})-(\ref{ex1_4}). Then we get:
\begin{equation}\label{ex1_5a} \frac{\partial \tilde n}{\partial v_1} = e^{v_1}\frac{\partial \tilde x}{\partial v_2},
\end{equation}
\begin{equation}\label{ex1_5b}
\frac{\partial \tilde n}{\partial v_2} = e^{-v_1} \frac{\partial \tilde x}{\partial v_1} +
\frac{e^{v_1}\bar b_{22}-e^{-v_1}\bar b_{11}}{\bar b_{12}} \frac{\partial \tilde x}{\partial v_2}.
\end{equation}
Viewing (\ref{ex1_5a})-(\ref{ex1_5b}) as Weingarten equations for
$\tilde F^2$, the shape operator of $\tilde F^2$ is given by
\begin{equation*}\label{ex1_6}
\tilde W = \left(\begin{array}{cc}0& e^{v_1} \\ e^{-v_1}& \frac{e^{v_1}\bar b_{22}-e^{-v_1}\bar
b_{11}}{\bar b_{12}}\end{array}\right).
\end{equation*}
The determinant of $\tilde W$, which is just the Gauss curvature of $\tilde F^2$, is equal to $-1$,
q.e.d.

\medskip

\textsl{Remark 4.2.} Notice that the Gauss curvature of the Riemannian metric (\ref{ex1_1}) is
equal to $- \frac{1}{(1-e^{-2v_1})^2}$ and therefore it is negative. The classical theory of
isometric immersions guarantees the existence of surfaces in $\mathbb R^3$ whose first fundamental
forms coincide with (\ref{ex1_1}). Besides, in the general case the assumption (\ref{ex1_1ss}) is
fulfilled too. Thus, two-dimensional surfaces in $\mathbb R^3$ satisfying both (\ref{ex1_1}) and
(\ref{ex1_1ss}) do exist and hence generate three-dimensional pseudo-spherical submanifolds in
$\mathbb R^5$ with degenerate Bianchi trans\-for\-ma\-ti\-ons of rank 2 as discussed above.

%\medskip

\textsl{Remark 4.3.} Notice that the horospherical coordinates
$(v_1, v_2, v_3)$ on the submanifold $F^3$ are subject to the
constraint $v_1>0$ which assures the positiveness of
(\ref{ex1_1}). Therefore,  $F^3$ represents a domain in a horoball
of the hyperbolic space isometrically immersed into $\mathbb R^5$.

\medskip

Let us calculate the fundamental forms of the submanifold $F^3\subset \mathbb R^5$ represented by
(\ref{ex1_2}). The vectors tangent to $F^3$ are \begin{equation*}\label{ex1_7a} \frac{\partial
x}{\partial v_1} = \left( \frac{\partial \bar x}{\partial v_1}, - e^{-v_1}\cos v_3, - e^{-v_1} \sin
v_3 \right),
\end{equation*}
\begin{equation*}\label{ex1_7b}
\frac{\partial x}{\partial v_2} = \left(\frac{\partial \bar x}{\partial v_2}, 0, 0 \right),
\end{equation*}
\begin{equation*}\label{ex1_7c}
\frac{\partial x}{\partial v_3} = \left(0, - e^{-v_1}\sin v_3,  e^{-v_1} \cos v_3 \right).
\end{equation*}
Taking into account (\ref{ex1_1}), it is easy to demonstrate that
the vectors
\begin{equation*}\label{ex1_8a}
n_1 = \left( \frac{e^{-v_1}}{\sqrt{1-e^{-v_1}}}\frac{\partial \bar
x}{\partial v_1}, \sqrt{1-e^{-v_1}}\cos v_3, \sqrt{1-e^{-v_1}}
\sin v_3 \right),
\end{equation*}
\begin{equation*}\label{ex1_8b}
n_2 = \left(\bar n, 0, 0 \right)
\end{equation*}
form an orthonormal frame in the normal plane of $F^3$.

Consequently, the second fundamental forms of $F^3$ with respect
to $n_1$ and $n_2$ read as follows:
\begin{equation*}\label{ex1_9a}
b^1 = \frac{e^{-v_1}}{\sqrt{1-e^{-2v_1}}} dv_1^2 +
\frac{e^{-3v_1}}{\sqrt{1-e^{-2v_1}}} dv_2^2 -
e^{-v_1}\sqrt{1-e^{-2v_1}} dv_3^2,
\end{equation*}
\begin{equation*}\label{ex1_9b}
b^2 =  \bar b_{11}dv_1^2 + 2\bar b_{12}dv_1dv_2 + \bar b_{22}
dv_2^2.
\end{equation*}
Moreover, the torsion coefficients of $F^3$, which are $\displaystyle{\mu_{12\vert j}=\langle
\frac{\partial n_1}{\partial v_j},n_2\rangle = - \langle n_1, \frac{\partial n_2}{\partial
v_j}\rangle}$ by definition, are expressed as follows
\begin{equation*}\label{ex1_10}
\mu_{12\vert 1}  = \frac{e^{-v_1}}{\sqrt{1-e^{-2v_1}}} \bar
b_{11}, \,\, \mu_{12\vert 2} = \frac{e^{-v_1}}{\sqrt{1-e^{-2v_1}}}
\bar b_{12} , \,\, \mu_{12\vert 3}  = 0.
\end{equation*}

The Gauss-Codazzi-Ricci equations for the submanifold $F^3\subset
\mathbb R^5$ reduce to the following system of equations:
\begin{equation*}\label{ex1_f1}
\bar b_{11}\bar b_{22} - \bar b_{12}^2 = -
\frac{e^{-2v_1}}{1-e^{-2v_1}},
\end{equation*}
\begin{equation*}\label{ex1_f2}
\frac{\partial \bar b_{11}}{\partial v_2} - \frac{\partial \bar
b_{12}}{\partial v_1} +  \frac{1}{1-e^{-2v_1}}\bar b_{12} = 0,
\end{equation*}
\begin{equation*}\label{ex1_f3}
\frac{\partial \bar b_{12}}{\partial v_2} - \frac{\partial \bar
b_{22}}{\partial v_1} -  \frac{e^{-2v_1}}{1-e^{-2v_1}}\bar b_{11}
- \bar b_{22} = 0.
\end{equation*}

These equations are written in terms of $\bar b_{11}$, $\bar
b_{12}$, $\bar b_{22}$, and they are just the Gauss-Codazzi
equations for the surface $\bar F^2\subset \mathbb R^3$.

\bigskip

\textbf{5. Example 2}

\medskip

In this section we will construct another kind of particular
three-dimensional pseudo-spherical submanifolds in $\mathbb R^5$
whose Bianchi trans\-for\-ma\-ti\-ons are degenerate of rank 2 and
lead to  two-dimensional surfaces.

Let $\bar F^2$ be a two-dimensional surface in $\mathbb R^3$
represented by a vector-function $\bar x(v_1,v_2)$ so that the
first fundamental form of $\bar F^2$ reads
\begin{equation}\label{ex2_1} d\bar x^2 = (1 - e^{-2v_1}\frac{v_2^2}{a^2})dv_1^2 + 2e^{-2v_1}\frac{v_2}{a^2}dv_1dv_2
+ e^{-2v_1}(1-\frac{1}{a^2})dv_2^2,
\end{equation}
where $a>1$ is an arbitrary constant.

Moreover, suppose the second fundamental form $\bar b = \bar
b_{11}dv_1^2 + 2\bar b_{12}dv_1dv_2 + \bar b_{22} dv_2^2$ of $\bar
F^2$  is not diagonalized,
\begin{equation}\label{ex2_1ss}
\bar b_{12}\not= 0.
\end{equation}.

The surface $\bar F^2\subset \mathbb R^3$ being given, consider a
three-dimensional submanifold $F^3$ in $\mathbb R^5$ represented
by the vector-function
\begin{equation}\label{ex2_2a}
x(v_1,v_2,v_3) = \left(\bar x(v_1,v_2), \frac{1}{a}e^{-v_1}v_2\cos
a v_3, \frac{1}{a} e^{-v_1} v_2\sin av_3 \right).
\end{equation}

Because of (\ref{ex2_1}), the first fundamental form of $F^3$
reads
\begin{equation}\label{ex2_2q}
dx^2 = dv_1^2 + e^{-2v_1}(dv_2^2 + v_2^2dv_3^2).
\end{equation}
If we introduce the coordinates $u_1=v_1$, $u_2=v_2\cos v_3$, $u_3=v_3\sin v_3$, then we get $dx^2
= du_1^2 + e^{-2u_1}(du_2^2 + du_3^2)$. Therefore, the submanifold $F^3$ is pseudo-spherical, its
Gauss curvature is equal to -1, and $(u_1, u_2, u_3)$ are horospherical coordinates on $F^3$. The
coordinates $(v_1, v_2, v_3)$ can be referred to as polar horospherical, since $(v_2, v_3)$
represent polar coordinates on each horosphere $v_1=const$, whereas $(u_2, u_3)$ are Cartesian
coordinates.

\medskip

{\bf Proposition 5.1.} {\it The Bianchi trans\-for\-ma\-ti\-on of
$F^3\subset \mathbb R^5$ corresponding to the choice of
horosherical coordinates is degenerate and transforms $F^3$ to a
two-dimensional surface $\tilde F^2\subset \mathbb R^5$ with the
following properties:

(i) $\tilde F^2$ belongs to a three-dimensional subspace $\mathbb
R^3$ of $\mathbb R^5$;

(ii) $\tilde F^2$ is pseudo-spherical, its Gauss curvature is equal to
$\displaystyle{-\frac{a^2}{a^2-1}}$. }

\medskip

\texttt{Proof.} Applying the Bianchi trans\-for\-ma\-ti\-on
(\ref{def}) to the vector-function (\ref{ex2_2a}),  it is easy to
get
\begin{equation}\label{ex2_2}
\tilde x = \left(\bar x+ \frac{\partial \bar x}{\partial v_1}, 0,
0 \right).
\end{equation}
Hence, the vector-function $\tilde x$ depends only on $v_1$, $v_2$ and takes values in the subspace
$\mathbb R^3$ of $\mathbb R^5$ determined by $x_4=0$, $x_5=0$.

Differentiate (\ref{ex2_2}) and apply the Weingarten equations
corresponding to $\bar F^2$. Then we have:
\begin{equation}\label{ex2_3a} \frac{\partial \tilde x}{\partial v_1} =  (1+\bar\Gamma^1_{11})
\frac{\partial \bar x}{\partial v_1} + \bar\Gamma^2_{11}
\frac{\partial \bar x}{\partial v_2} + \bar b_{11} \bar n =
\frac{a^2-1}{a^2-1 - v_2^2e^{-2v_1}}\frac{\partial \bar
x}{\partial v_1} - \frac{v_2}{a^2-1 -
v_2^2e^{-2v_1}}\frac{\partial \bar x}{\partial v_2} + \bar b_{11}
\bar n,
\end{equation}
\begin{equation}\label{ex2_3b}
\frac{\partial \tilde x}{\partial v_2} = \bar\Gamma^1_{12}
\frac{\partial \bar x}{\partial v_1} + (1+\bar\Gamma^2_{12})
\frac{\partial \bar x}{\partial v_2} + \bar b_{12} \bar n = \bar
b_{12} \bar n,
\end{equation}
where $\bar\Gamma^1_{11}=\frac{v_2^2e^{-2v_1}}{a^2-1 - v_2^2e^{-2v_1}}$,
$\bar\Gamma^2_{11}=-\frac{v_2}{a^2-1 - v_2^2e^{-2v_1}}$, $\bar\Gamma^1_{12}=0$,
$\bar\Gamma^2_{12}=-1$ are the Christoffel symbols and $\bar n$ is the unit normal of $\bar F^2$.
Taking into account (\ref{ex2_1ss}), we see that $\displaystyle{\frac{\partial \tilde x}{\partial
v_2}}$ and $\displaystyle{\frac{\partial \tilde x}{\partial v_2}}$ are non-collinear.

Thus, $\tilde x(v_1,v_2)$ given by (\ref{ex2_2}) represents a
two-dimensional surface, $\tilde F^2$, in the three-dimensional
Euclidean space $\mathbb R^3$.

To finish the proof, let us find the Gauss curvature $\tilde K$ of $\tilde F^2$. It is easy to
verify in view of (\ref{ex2_1}), (\ref{ex2_3a})-(\ref{ex2_3b}) that the unit normal to $\tilde
F^2\subset \mathbb R^3$ is written as follows:
\begin{equation}\label{ex2_4}
\tilde n = \frac{e^{v_1}}{\sqrt{1-\frac{1}{a^2}}} \frac{\partial
\bar x}{\partial v_2}.
\end{equation}
Differentiate (\ref{ex2_4}), apply the Weingarten equations of $\bar F^2$ and express
$\displaystyle{\frac{\partial \bar x}{\partial v_1}}$, $\displaystyle{\frac{\partial \bar
x}{\partial v_2}}$, $\bar n$ in terms of $\displaystyle{\frac{\partial \tilde x}{\partial v_1}}$,
$\displaystyle{\frac{\partial \tilde x}{\partial v_2}}$, $\tilde n$ by using
(\ref{ex2_3a})-(\ref{ex2_4}). Then we get:
\begin{equation}\label{ex2_5a}
\frac{\partial \tilde n}{\partial v_1} =
\frac{e^{v_1}}{\sqrt{1-\frac{1}{a^2}}}\frac{\partial \tilde
x}{\partial v_2},
\end{equation}
\begin{equation}\label{ex2_5b}
\frac{\partial \tilde n}{\partial v_2} =
\frac{e^{-v_1}}{\sqrt{1-\frac{1}{a^2}}}\frac{\partial \tilde
x}{\partial v_1} + \frac{e^{v_1}\bar b_{22}-e^{-v_1}\bar
b_{11}}{\bar b_{12}\sqrt{1-\frac{1}{a^2}}} \frac{\partial \tilde
x}{\partial v_2}.
\end{equation}
Viewing (\ref{ex2_5a})-(\ref{ex2_5b}) as Weingarten equations for
$\tilde F^2$, the shape operator of $\tilde F^2$ is given by
\begin{equation*}\label{ex2_6}
\tilde W = \left(\begin{array}{cc}0&
\frac{e^{v_1}}{\sqrt{1-\frac{1}{a^2}}} \\
\frac{e^{-v_1}}{\sqrt{1-\frac{1}{a^2}}}& \frac{e^{v_1}\bar
b_{22}-e^{-v_1}\bar b_{11}}{\bar
b_{12}\sqrt{1-\frac{1}{a^2}}}\end{array}\right).
\end{equation*}
The determinant of $\tilde W$, which is just the Gauss curvature of $\tilde F^2$, is equal to
$\displaystyle{-\frac{a^2}{a^2-1}}$, q.e.d.

\medskip

\textsl{Remark 5.2.} Notice that the Gauss curvature of the Riemannian metric (\ref{ex2_1}) is
equal to

\noindent $-\frac{a^2(a^2-1)}{(a^2-1-v_2^2e^{-2v_1})^2}$ and therefore it is negative. The
classical theory of isometric immersions guarantees the existence of surfaces in $\mathbb R^3$
whose first fundamental forms coincide with (\ref{ex2_1}). Besides, in the general case the
assumption (\ref{ex2_1ss}) is fulfilled too. Thus, two-dimensional surfaces in $\mathbb R^3$
satisfying both (\ref{ex2_1}) and (\ref{ex2_1ss}) do exist and hence generate three-dimensional
pseudo-spherical submanifolds in $\mathbb R^5$ with degenerate Bianchi trans\-for\-ma\-ti\-ons of
rank 2 as discussed above.

\medskip

\textsl{Remark 5.3.} The polar horospherical coordinates
$(v_1,v_2,v_3)$ are subject to the constraint $v_2
e^{-v_1}<\sqrt{a^2-1}$, which assures the positiveness of
(\ref{ex2_1}). Therefore, $F^3$ in question represents a specific
cone-like domain of the hyperbolic space isometrically immersed
into $\mathbb R^5$.

\medskip

Let us calculate the fundamental forms of the submanifold
$F^3\subset \mathbb R^5$ represented by (\ref{ex2_2}). The vectors
tangent to $F^3$ are
\begin{equation*}\label{ex2_7a} \frac{\partial
x}{\partial v_1} = \left( \frac{\partial \bar x}{\partial v_1}, -
\frac{1}{a} e^{-v_1} v_2 \cos a v_3, - \frac{1}{a} e^{-v_1} v_2
\sin a v_3 \right),
\end{equation*}
\begin{equation*}\label{ex2_7b}
\frac{\partial x}{\partial v_2} = \left(\frac{\partial \bar
x}{\partial v_2}, \frac{1}{a} e^{-v_1} \cos a v_3, \frac{1}{a}
e^{-v_1} \sin a v_3 \right),
\end{equation*}
\begin{equation*}\label{ex2_7c}
\frac{\partial x}{\partial v_3} = \left(0, - e^{-v_1} v_2 \sin a
v_3, e^{-v_1} v_2 \cos a v_3 \right).
\end{equation*}
Taking into account (\ref{ex2_1}), it is easy to demonstrate that
the vectors
\begin{equation*}\label{ex2_8a}
n_1 = \frac{\sqrt{a^2-1-e^{-2v_1} v_2^2}}{a} \left(
\frac{a}{{a^2-1-e^{-2v_1} v_2^2}}( e^{-v_1}v_2\frac{\partial \bar
x}{\partial v_1} - e^{v_1}\frac{\partial \bar x}{\partial v_2}),
\cos a v_3, \sin av_3 \right),
\end{equation*}
\begin{equation*}\label{ex2_8b}
n_2 = \left(\bar n, 0, 0 \right)
\end{equation*}
form an orthonormal frame in the normal plane of $F^3$.

Consequently, the second fundamental forms of $F^3$ with respect
to $n_1$ and $n_2$ read as follows:
\begin{equation*}\label{ex2_9a}
b^1 = \frac{e^{-v_1}v_2}{\sqrt{a^2-1-e^{-2v_1}v_2^2}}\, dv_1^2 +
\frac{e^{-3v_1}v_2}{\sqrt{a^2-1-e^{-2v_1}v_2^2}}\, dv_2^2 -
e^{-v_1}v_2\sqrt{a^2-1-e^{-2v_1}v_2^2}\, dv_3^2,
\end{equation*}
\begin{equation*}\label{ex2_9b}
b^2 =  \bar b_{11}\, dv_1^2 + 2\bar b_{12}\, dv_1dv_2 + \bar b_{22}\, dv_2^2.
\end{equation*}
Moreover, the torsion coefficients of $F^3$ are expressed as follows:
\begin{equation*}\label{ex2_10}
\mu_{12\vert 1}  = \frac{e^{-v_1}v_2\bar b_{11}-e^{v_1}\bar
b_{12}}{\sqrt{a^2-1-e^{-2v_1}v_2^2}} \, ,\,\, \mu_{12\vert 2} =
\frac{e^{-v_1}v_2\bar b_{12}-e^{v_1}\bar
b_{22}}{\sqrt{a^2-1-e^{-2v_1}v_2^2}}\, , \,\,\mu_{12\vert 3} = 0.
\end{equation*}

The Gauss-Codazzi-Ricci equations for the submanifold $F^3\subset
\mathbb R^5$ reduce to the following system of equations:
\begin{equation*}\label{ex2_f1}
\bar b_{11}\bar b_{22} - \bar b_{12}^2 = -
\frac{(a^2-1)e^{-2v_1}}{a^2-1-v_2^2e^{-2v_1}},
\end{equation*}
\begin{equation*}\label{ex2_f2}
\frac{\partial \bar b_{11}}{\partial v_2} - \frac{\partial \bar
b_{12}}{\partial v_1} +
\frac{1}{a^2-1-v_2^2e^{-2v_1}}\left((a^2-1)\bar b_{12}-v_2\bar
b_{22}\right) = 0,
\end{equation*}
\begin{equation*}\label{ex2_f3}
\frac{\partial \bar b_{12}}{\partial v_2} - \frac{\partial \bar
b_{22}}{\partial v_1} -
\frac{e^{-2v_1}}{a^2-1-v_2^2e^{-2v_1}}\left( (a^2-1)\bar
b_{11}-v_2\bar b_{12}\right)  - \bar b_{22} = 0.
\end{equation*}

These equations are written in terms of $\bar b_{11}$, $\bar
b_{12}$, $\bar b_{22}$, and they are just the Gauss-Codazzi
equations for the surface $\bar F^2\subset \mathbb R^3$.

\medskip

\textsl{Remark 5.4.} Since the parameter $a\in (1,+\infty)$ can be fixed arbitrarily, we
constructed not a single submanifold but a one-parameter family of specific three-dimensional
pseudo-spherical submanifolds in $\mathbb R^5$ with degenerate Bianchi trans\-for\-ma\-ti\-ons of
rank 2. The dependence on $a$ is essential and can not be destroyed by scalings of coordinates,
because the Gauss curvature $\tilde K$ of the resulting surfaces in $\mathbb R^3$ is a strongly
increasing function of $a$. Notice also that $\tilde K$ is less than $-1$, and it tends to $-1$ as
$a\to +\infty$.

\bigskip

\textbf{6. Proof of Theorem}

\medskip

To complete the proof of our main Theorem, we will use the
uniqueness part of the generalized fundamental theorem (Bonnet
theorem) which claims the following: if two $n$-dimensional
submanifolds in the $(n+m)$-dimensional Euclidean space have the
same first fundamental form, second fundamental forms and torsion
coefficients, then these submanifolds coincide up to a rigid
motion in the ambient space, c.f. \cite{T1}.

In Sections 2-3 it was demonstrated that if a pseudo-spherical
submanifold $F^3$ in $\mathbb R^5$ admits a
ho\-lo\-no\-mi\-cal\-ly degenerate Bianchi trans\-for\-ma\-ti\-on
of rank 2, then it can be parameterized by specific coordinates
($v_1$, $v_2$, $v_3$) and equipped with a specific normal frame
$n_1$, $n_2$ so that either its fundamental forms are
\begin{equation*}
ds^2 \, =\, dv_1^2 \, +\, e^{-2v_1} \left(dv_2^2+dv_3^2\right),
\end{equation*}
\begin{equation*}
{\text b^1} \, =\, \frac{1}{\sqrt{e^{2v_1}-1}}\, dv_1^2 + \frac{e^{-2v_1}}{\sqrt{e^{2v_1}-1}} \,
dv_2^2 -e^{-2v_1}\sqrt{e^{2v_1}-1}\, dv_3^2, \end{equation*}
\begin{equation*}
{\text b^2}  \, =\,  {\text b^2_{11}}\, dv_1^2 + 2{\text b^2_{12}}\, dv_1dv_2 + {\text b^2_{22}}\,
dv_2^2,
\end{equation*}
\begin{equation*}
\mu_{12\vert 1}  \, =\, \frac{{\text b^2_{11}}} {\sqrt{e^{2v_1}-1}},\, \mu_{12\vert 2} \, =\,
\frac{{\text b^2_{12}}} {\sqrt{e^{2v_1}-1}},\, \mu_{12\vert 3}  \, =\,  0,
\end{equation*}
or its fundamental forms are
\begin{equation*}
ds^2  \, =\,  dv_1^2+e^{-2v_1} \left(dv_2^2+v_2^2dv_3^2\right),
\end{equation*}
\begin{equation*}
{\text b^1}  \, =\,  \frac{v_2}{\sqrt{f_0e^{2v_1}-(v_2)^2}}\, dv_1^2 +
e^{-2v_1}\frac{v_2}{\sqrt{f_0e^{2v_1}-(v_2)^2}}\, dv_2^2 - e^{-2v_1}v_2\sqrt{f_0e^{2v_1}-(v_2)^2}\,
dv_3^2,
\end{equation*}
\begin{equation*}
{\text b^2}  \, =\,  {\text b^2_{11}}\, dv_1^2 + 2{\text b^2_{12}}\, dv_1dv_2 + {\text b^2_{22}}\,
dv_2^2,
\end{equation*}
\begin{equation*}
\mu_{12\vert 1}  \, =\, \frac{v_2 {\text b^2_{11}} - e^{2v_1}{\text
b^2_{12}}}{\sqrt{f_0e^{2v_1}-(v_2)^2}}, \, \mu_{12\vert 2}  \, =\,  \frac{v_2 {\text b^2_{12}} -
e^{2v_1}{\text b^2_{22}}}{\sqrt{f_0e^{2v_1}-(v_2)^2}}, \, \mu_{12\vert 3}  \, =\,  0.
\end{equation*}

In the first case $F^3$ has the same fundamental forms as the
submanifold described in Section 4. By the Bonnet fundamental
theorem mentioned above, since the fundamental forms are the same,
then the submanifolds coincide up to a rigid motion in $\mathbb
R^5$. Hence, $F^3$ has the properties described in Proposition
4.1, i.e., the two-dimensional surface obtained from $F^3$ by the
degenerate Bianchi transformation belongs to a subspace $\mathbb
R^3\subset \mathbb R^5$ and its Gauss curvature is equal to $-1$.

Similarly, in the second case $F^3$ has the same fundamental forms as the submanifold described in
Section 5, if one sets $f_0=a^2-1$. By the Bonnet fundamental theorem, since the fundamental forms
are the same, then the submanifolds coincide up to a rigid motion in $\mathbb R^5$. Hence, $F^3$
has the properties described in Proposition 5.1, i.e., the two-dimensional surface obtained from
$F^3$ by the degenerate Bianchi transformation belongs to a subspace $\mathbb R^3\subset \mathbb
R^5$ and has constant negative Gauss curvature less than $-1$. This completes the proof of Theorem.

\end{document}